\documentclass[leqno,12pt]{amsart}
\usepackage{amsfonts}
\usepackage{amsmath,amssymb,amsthm}
\usepackage{amsfonts}

\newtheorem{theorem}{Theorem}[section]
\newtheorem{proposition}[theorem]{Proposition}
\newtheorem{lemma}[theorem]{Lemma}
\newtheorem{corry}[theorem]{Corollary}
\newtheorem{defi}[theorem]{Definition}
\theoremstyle{definition}
\newtheorem{example}[theorem]{Example}
\theoremstyle{definition}
\newtheorem{rem}[theorem]{Remark}

\numberwithin{claim}{theorem}

\renewenvironment{proof}{\textit{Proof.}}{\hfill\ensuremath{\qed}}

\def \qed{\hfill{\hbox{$\square$}}}

\numberwithin{equation}{section}

\begin{document}
\title[Class $\mathcal A$ Surfaces in $L^4_1(f,0)$]{On  Space-like Class $\mathcal A$ Surfaces in Robertson-Walker Space Times}

\author[B. Bekta\c{s} Dem\.{i}rc\.i]{Burcu Bekta\c s Dem\.{i}rc\.i}
\address{Fatih Sultan Mehmet Vak{\i}f University, Hal\.{I}\c{c} Campus, Faculty of Engineering,
Department of Software Engineering, 34445, Beyo\u{g}lu, \.{I}stanbul, T{\"u}rk\.{I}ye}
\email{bbektas@fsm.edu.tr, 0000-0002-5611-5478}

\author[N. Cenk Turgay]{Nurettin Cenk Turgay}
\address{Department of Mathematics, Faculty of Science and Letters, Istanbul Technical University, \.{I}stanbul, T{\"u}rk\.{I}ye}
\email{nturgay@itu.edu.tr, 0000-0002-0171-3876}

\author[R. Ye\u{g}\.{i}n \c{S}en]{R\"{u}ya Ye\u{g}\.{i}n \c{S}en}
\address{Department of Mathematics, Faculty of Engineering and Natural Sciences, \.{I}stanbul Medeniyet University, \.{I}stanbul, T{\"u}rk\.{I}ye}
\email{ruya.yegin@medeniyet.edu.tr, 0000-0002-2642-1722}

\maketitle

\begin{abstract}
In this article, we consider space-like surfaces in Robertson-Walker Space times $L^4_1(f,c)$ with comoving observer field $\frac{\partial}{\partial t}$.
We study some problems related to such surfaces satisfying  
the geometric conditions imposed on the tangential part and normal part of the unit vector field 
$\frac{\partial}{\partial t}$ naturally defined. 
First, we investigate space-like surfaces in $L^4_1(f,c)$ satisfying that 
the tangent component of $\frac{\partial}{\partial t}$ is an eigenvector of all shape operators, called class $\mathcal A$ surfaces. 
Then, we get a classification theorem of space-like class $\mathcal A$ surfaces in $L^4_1(f,0)$. 
Also, we examine minimal space-like class $\mathcal A$ surfaces in $L^4_1(f,0)$. 
Finally, we give the parametrizations of space-like surfaces in $L^4_1(f,0)$ 
when the normal part of the unit vector field $\frac{\partial}{\partial t}$ is parallel. 
\end{abstract}

\section{Introduction}
In recent years, there has been significant interest among geometers in studying submanifolds of product spaces 
resulting in numerous findings. Some of them are given by \cite{Fetcu15}, \cite{Anciaux15}, \cite{ChenWei2009}, \cite{Dillen12}, \cite{Dillen09}.

Apart from Cartesian product spaces, the other example can be considered as Robertson Walker-Space times which are 4-dimensional Lorentzian manifolds. 
The Robertson Walker-Space times, denoted by $L^4_1(f,c)$, are defined as Cartesian products of space forms by a real interval equipped 
with a Lorentzian warped product metric. 
Thus, the family of Robertson Walker space times includes the de Sitter, Minkowski, and the anti-de Sitter space time
and also Friedmann's cosmological models. 
In physics, Robertson Walker-Space times are important due to the fact that they explain homogeneous,
isotropic expanding and contracting universes, see \cite{Hawking1973} and \cite{ONeill1982}.

From geometrical point of view, there are some studies related to classification of surfaces or hypersurfaces in Robertson Walker space times, 
see \cite{Alias95}, \cite{Alias96}, \cite{Alias97}, \cite{AnciauxCipriani2020}, \cite{DekimpeJVanderVeken2020}.
Especially, B.-Y. Chen and J. Van der Veken \cite{ChenJVanderVeken2007} investigated space-like surfaces and Lorentzian surfaces 
in Robertson Walker space times having some important geometric properties such 
as marginally trapped, positive relative nullity and totally geodesic, etc. 

In product spaces, there exists a unit vector field spanning the factor $\mathbb{R}$, denoted by  $\frac{\partial}{\partial t}$, 
and B. Mendon\c{c}a and R. Tojeiro \cite{MendonTojeiro2014}  mentioned the decomposition of it given as follows.

Given an isometric immersion $\phi:M\rightarrow\mathbb{Q}_{c}^n\times\mathbb{R}$, 
a tangent vector field $T$ on $M$ and a normal vector field $\eta$ along $\phi$ are defined by
\begin{equation}\label{RWTS-etaandTDef} 
\frac{\partial}{\partial t}=\phi_* T+\eta.    
\end{equation}
Here, $\mathbb{Q}_c^n$ denotes the space forms.
For the Robertson walker space times, there still exists the unit vector field $\frac{\partial}{\partial t}$ which is tangent to first factor and 
it can be decomposed in a similar way as \eqref{RWTS-etaandTDef}. Also, this unit vector field is known as a comoving observer field in general relativity, \cite{ChenJVanderVeken2007}.

In \cite{Tojeiro2010}, R. Tojeiro studied hypersurfaces in product spaces $\mathbb{Q}_{c}^n\times\mathbb{R}$
for which the tangent component of $\frac{\partial}{\partial t}$ in \eqref{RWTS-etaandTDef} is an eigenvector of all shape operators.
Moreover, B. Mendon\c{c}a and R. Tojeiro \cite{MendonTojeiro2014} obtained the characterization of
all such submanifolds in the product spaces $\mathbb{Q}_{c}^n\times\mathbb{R}$.
In this context, the following definition based on \cite{Tojeiro2010} and \cite{MendonTojeiro2014} for submanifolds of the Cartesian product space
can be given. 
\begin{defi}
\label{defClassA}
If the vector $T$ in \eqref{RWTS-etaandTDef} is an eigenvector of all shape operators of $M$, 
then $M$ is said to be a class $\mathcal A$ surface. 
\end{defi}

In this article, we study space-like surfaces in Robertson-Walker space times satisfying certain properties 
in terms of the vectors $T$ and $\eta$ in the decomposition \eqref{RWTS-etaandTDef}. 
First, we investigate space-like class $\mathcal A$ surfaces in $L^4_1(f,c)$ and then, 
we give a local classification theorem for space-like class $\mathcal A$ surfaces in $L^4_1(f,0)$.
Also, we determine space-like class $\mathcal A$ surfaces in $L^4_1(f,0)$ with zero mean curvature vector. 
Finally, we prove that the vector field $\eta$ in \eqref{RWTS-etaandTDef} is parallel if and only if the space-like surface in $L^4_1(f,c)$ must be an 
element of class $\mathcal A$ surfaces. By using this, we obtain the parametrizations of space-like surfaces in $L^4_1(f,0)$ 
with parallel vector field $\eta$.


\section {Preliminaries}

Let ${\mathbb{Q}}^{n-1}_c$ denote the $n-1$ dimensional Riemannian space-form with the constant sectional curvature $c$, i.e.,
$${\mathbb{Q}}^{n-1}_c=\left\{\begin{array}{cc}
\mathbb S^{n-1}&\mbox{if }c=1,\\
\mathbb E^{n-1}&\mbox{if }c=0,\\
\mathbb H^{n-1}&\mbox{if }c=-1
\end{array}\right.$$
and $g_c$ stand for its metric tensor.

If $I$ is an open interval and $f:I\to\mathbb R$ is a smooth, non-vanishing function, 
then the Robertson-Walker space time $L^4_1(f,c)$ is defined as the Lorentzian warped product $I^1_1\times_f {\mathbb {Q}}^{3}_c$ 
whose metric tensor $\tilde g$ is
$$\tilde g=-dt^2+f^2(t)g_c.$$

Let  $\Pi_1:I\times \mathbb{Q}^{3}_c\to I$ and $\Pi_2:I\times \mathbb{Q}^{3}_c\to \mathbb{Q}^{3}_c$
denote the canonical projections. 
For a given vector field $X$ in ${L}^4_1(f,c)$, 
we define a function $X_0$ and a vector field $\bar X$ by the decomposition
$$X=X_0\frac{\partial}{\partial t}+\bar X.$$
Note that we have $X_0=-\tilde g\left(\frac{\partial}{\partial t},X \right)$ and $\Pi_1^*(\bar X)=0$.

First, we would like to express the Levi-Civita connection of $ L^4_1(f,c)$. 
Note that the following lemma can be directly obtained from \cite[Lemma 2.1]{ChenJVanderVeken2007}.
\begin{lemma}\label{LemmaLn1f0LCConnect}
The Levi-Civita connection $\widetilde\nabla$ of ${L}^4_1(f,c)$ is
\begin{eqnarray}
\label{RWSTconn}
\widetilde{\nabla}_X Y&=&\nabla^0_XY+ (\ln{f})' \left(\tilde g (\bar X,\bar Y){\partial_t}+X_0\bar Y+Y_0\bar X\right)
\end{eqnarray}
whenever $X$ and $Y$ are tangent to ${L}^4_1(f,c)$, 
where $\nabla^0$ denotes the Levi-Civita connection of the Cartesian product space ${L}^4_1(1,c)=I\times \mathbb{Q}^{3}_c$.
\end{lemma}

\subsection{Space-like surfaces in $ L^4_1(f,c)$}
Consider an oriented space-like surface  $M$ in $L^4_1(f,c)$ with the Levi-Civita connection $\nabla$ and metric tensor $g$. 
Through the misuse of terminology, we shall denote the induced connection of $L^4_1(f,c)$ by $\widetilde\nabla$. Then, the Gauss and 
Weingarten formulas 
\begin{eqnarray}
	\label{Gauss} \widetilde\nabla_X Y&=& \nabla_X Y + h(X,Y),\\
	\label{Weingarten} \widetilde\nabla_X \xi&=& -A_\xi(X) +\nabla^\perp_X \xi
\end{eqnarray}
define  the second 
fundamental form  $h$, the shape operator $A$  and the normal connection $\nabla^\perp$ of $M$, where $X,Y$ are tangent to $M$ and $\xi$ is normal to $M$. 
Note that $A$ and $h$ are related by
\begin{eqnarray}
	\label{AhRelatedBy} 
 \tilde g(h(X,Y),\xi)=g(A_\xi X,Y).
\end{eqnarray}

Let $\phi:\Omega\to  L^4_1(f,c)$ be a local parametrization of $M$ and put $\mathcal T=\Pi_1\circ\phi$. 
\begin{rem}
\label{remT}
If $\mathrm{grad\,} \mathcal{T}=0$ on $M$, then we have $(M,g)\subset \{ t_0\} \times_{f(t_0)}{\mathbb{Q}}^3_c$. 
Throughout this work, we are going to exclude this case and assume the existence of $p\in M$ at which $\mathrm{grad\,} \mathcal{T}\neq0$.
\end{rem}

Since $M$ is a space-like surface, by considering the decomposition \eqref{RWTS-etaandTDef} one may define a function $\theta$ and a positively oriented \textit{global} orthonormal frame field $\{e_1,e_2;e_3,e_4\}$ on $M$ by 
\begin{equation}\label{RWTS-etaandTDefNew}
\left.\frac{\partial}{\partial t}\right|_M=\sinh\theta \, e _1+\cosh\theta \, e _3.  
\end{equation}
where $g(e_1,e_1)=g(e_2, e_2)=1$ and $\tilde{g}(e_3,e_3)=-\tilde{g}(e_4,e_4)=-1$.
Throughout this article, we consider such an orthonormal frame field $\{e_1,e_2;e_3,e_4\}$ on $M$. 

Note that we are going to use the  notation 
$$h^\alpha_{ij}=\tilde g(h(e_i,e_j),e_\alpha)=g(e_i,A_{e_\alpha}e_j)$$
for the coefficients of the second fundamental form and
$\omega_{12}$ will stand for the connection form defined by
$$\omega_{12}(X)=g(\nabla_{X}e_1,e_2)=- g(\nabla_{X}e_2,e_1).$$
Then, the mean curvature vector of $M$ is defined by
$$H=\frac 12\mathrm{tr\,} h=-\frac{h^3_{11}+h^3_{22}}{2}e_3+ \frac{h^4_{11}+h^4_{22}}{2}e_4$$
and $M$ is said to be minimal if $H$ vanishes on $M$.
Also, we are going to use the following lemma:
\begin{lemma}\label{Lemmafptsubman}
Let $\theta$ and $e_1$ be as defined above. Then,
\begin{equation}\label{LemmafptsubmanEqu}
    \left.f'\right|_M=-\frac{1}{\sinh\theta}e_1\left(\left.f\right|_M\right).
\end{equation}
\end{lemma}
\begin{proof}
Let $p=(t_0,\tilde p)\in M$ and consider an integral curve $\alpha=(\alpha_0,\tilde \alpha)$ of $e_1$ starting from $p$. Then we  have
$$e_1(\left.f\right|_M)_p=\left.\frac{d}{du}\right|_{u=0} (f\circ\alpha)(u)=\left.\frac{d}{du}\right|_{u=0} f(\alpha_0(u))=\alpha'_0(0)f'(\alpha_0(0))$$
Since $\alpha(0)=p$ and $\alpha'(0)=(e_1)_p$, the last equation implies
$$e_1(\left.f\right|_M)_p=-\sinh(\theta(p))f'(t_0)$$
which yields \eqref{LemmafptsubmanEqu}.
\end{proof}

During this work, the manifolds that we are dealing with are smooth and connected unless otherwise is stated.

\section{Basic Facts for Space-like Surfaces in Robertson-Walker Space Times $L^4_1(f,c)$}
In this section, we give some basic facts about space-like surfaces in $L^4_1(f,c)$. 
\begin{lemma}
\label{LemmaCurvProps}
Let $M$ be a space-like surface in $L^4_1(f,c)$. Then, we have the followings:
\begin{subequations}\label{LemmaCurvPropsALL}
\begin{eqnarray}
\label{LemmaCurvProps1}\nabla_{e_1}e_1&=&h^3_{12}\coth\theta \, e _2,\\
\label{LemmaCurvProps2}\nabla_{e_2}e_1&=&\left((\ln {f})' \mathrm{csch\,}\theta+h^3_{22}\coth\theta \right)e_2,\\
\label{LemmaCurvProps5}\nabla^\perp_{e_1}e_3&=&-h^4_{11}\tanh\theta \, e_4,\\
\label{LemmaCurvProps6}\nabla^\perp_{e_2}e_3&=&-h^4_{12}\tanh\theta \, e_4\\
\label{LemmaCurvProps3}e_1(\theta)&=&(\ln {f}) '\cosh\theta+h^3_{11},\\
\label{LemmaCurvProps4}e_2(\theta)&=&h^3_{12}.
\end{eqnarray}
\end{subequations}
\end{lemma}

\begin{proof}
Assume that $M$ is a space-like surface in $L^4_1(f,c)$.
From \eqref{RWTS-etaandTDefNew}, we get $(e_2)_0=\tilde{g}\left(e_2,\frac{\partial}{\partial t}\right)=0$ 
and $(e_1)_0=-\tilde g\left(\frac{\partial}{\partial t},e_1\right)=-\sinh\theta$. 
Then, we have $e_2=\bar e_2$ and $e_1=-\sinh{\theta}\frac{\partial}{\partial t}+\bar{e_1}$.
Considering \eqref{RWSTconn} in Lemma \ref{LemmaLn1f0LCConnect} with these equations, one can obtain
\begin{eqnarray}
\label{TNe1pd}\widetilde\nabla_{e_1} \frac{\partial}{\partial t}&=&  (\ln{f})' \left(\cosh^2\theta \, e _1+\sinh\theta\cosh\theta \, e_3\right),\\
\label{TNe2pd}\widetilde\nabla_{e_2} \frac{\partial}{\partial t}&=&  (\ln{f})' e_2.
\end{eqnarray}
On the other hand, the equation \eqref{RWTS-etaandTDefNew} implies that 
\begin{eqnarray*}
\widetilde\nabla_{X}\frac{\partial}{\partial t}= X(\theta) \cosh{\theta} e_1 + \sinh{\theta} \tilde{\nabla}_X e_1
+ X(\theta) \sinh{\theta} e_3 + \cosh{\theta} \tilde{\nabla}_X e_3
\end{eqnarray*}
from which we obtain
\begin{align}
\label{TNXpdTang}
\begin{split}\left(\widetilde\nabla_{X}\frac{\partial}{\partial t}\right)^T=&X(\theta)\cosh\theta e_1+\sinh\theta \nabla_{X}e _1-\cosh\theta A_{e_3}X,\\
\left(\widetilde\nabla_{X}\frac{\partial}{\partial t}\right)^\perp=&X(\theta)\sinh\theta e_3+\sinh\theta h(X,e_1)+\cosh\theta \nabla^\perp_X e_3.
\end{split}
\end{align}
Comparing the equation \eqref{TNXpdTang} for $X=e_1$ and \eqref{TNe1pd}, we get \eqref{LemmaCurvProps1}, \eqref{LemmaCurvProps5} and \eqref{LemmaCurvProps3}. Similarly, the equation \eqref{TNXpdTang} for $X=e_2$
and \eqref{TNe2pd} give \eqref{LemmaCurvProps2}, \eqref{LemmaCurvProps6} and \eqref{LemmaCurvProps4}.
\end{proof}

From Lemma \ref{LemmaCurvProps}, we give the following proposition. 
\begin{proposition}
\label{PropClassA}
Let $M$ be a space-like surface in $L^4_1(f,c)$. Then, the following statements are equivalent to each other:
\begin{enumerate}
\item [(i)] $M$ is space-like class $\mathcal A$ surface in $L^4_1(f,c)$.
\item [(ii)] For the vector fields $e_1, e_2$ tangent to $M$, we have $h(e_1,e_2)=0$. 
\item [(iii)] For the vector vector fields $e_2, e_4$, we have $e_2(\theta)=0$ and $\nabla^\perp_{e_2}e_4=0$. 
\end{enumerate}
\end{proposition}

\begin{proof}
From Definition \ref{defClassA} and \eqref{AhRelatedBy}, it can be easily seen that statements (i) and (ii) 
are equivalent to each other.
Statements (ii) and (iii) are direct consequences of \eqref{LemmaCurvProps6} and \eqref{LemmaCurvProps4}.
\end{proof}

\begin{lemma}
\label{space-likeSurfCanonicalParamLemma2}
Let $M$ be a space-like class $\mathcal A$ surface in $L^4_1(f,c)$ and $p\in M$.  
Then, there exists a local coordinate system $(u,v)$ defined on a neighborhood $\mathcal N_p$ of $p$ 
which can be parameterized by 
\begin{equation}\label{space-likeSurfCanonicalParam}
\phi(u,v)=(u,\tilde\phi(u,v)) 
\end{equation}
for an immersion $\tilde\phi:\Omega\subset\mathbb R^2\to \mathbb{Q}^{3}_c$ satisfying 
\begin{equation}
\label{space-likeSurfCanonicalParamLemma2Eq1}
g_c(\tilde\phi_u,\tilde\phi_v)=0\;\;\; \mbox{and}\;\; 
\partial_v (g_c(\tilde\phi_{u},\tilde\phi_u))=0.
\end{equation}
\end{lemma}
\begin{proof}
Assume that $M$ is a space-like class $\mathcal A$ surface in $L^4_1(f,c)$. 
Then, the equations \eqref{LemmaCurvProps1} and \eqref{LemmaCurvProps2} become $\nabla_{e_1}e_1=0$ 
and $\nabla_{e_2}e_1=(\ln {f})'\mathrm{csch\,}\theta e_2$. 
Define $X=-\frac 1{\sinh\theta}e_1$ and $Y=\gamma e_2$. 
Then, we have 
\begin{align}
[X,Y]=\mathrm{csch\,}\theta \left(-e_1(\gamma)+\gamma (\ln {f})'\mathrm{csch\,}\theta \right) e_2.
\end{align}
Thus, we choose a smooth function $\gamma$ satisfying $e_1(\gamma)=\gamma\omega$ which makes $[X,Y]=0.$
Hence, there exists a neighborhood $\mathcal N_p$ of $p$ on which a local coordinate system $(u,v)$ is defined such that
\begin{equation}
\label{space-likeSurfCanonicalParamLemma2Eq2}
X=-\frac 1{\sinh\theta}e_1=\partial_u,\qquad Y=\gamma e_2=\partial_v.
\end{equation}
Now, we consider a parametrization of $M$ in $L^4_1(f,c)$ 
\begin{equation}
    \phi(u,v)=(\mathcal T(u,v), \tilde{\phi}(u,v))
\end{equation}
where $\mbox{grad}{\mathcal T}\neq 0$ and $\tilde{\phi}:\Omega\subset\mathbb R^2\to \mathbb{Q}^{3}_c$ is an immersion. 
From the equations $\tilde g\left(e_1,\frac{\partial}{\partial t}\right)=\sinh\theta$ 
and $\tilde g\left(e_2,\frac{\partial}{\partial t}\right)=0$,
we have 
$$\frac{\partial \mathcal T}{\partial u}=1,\qquad \frac{\partial \mathcal T}{\partial v}=0$$
on $\mathcal N_p$, respectively. Thus, $M$ has a parametrization given by \eqref{space-likeSurfCanonicalParam}
with $g_c(\tilde\phi_u,\tilde\phi_v)=0$.
Moreover, we get
\begin{equation}
\label{lemmaeq1}
\tilde{g}(\phi_u,\phi_u)=1+f^2g_c(\tilde{\phi}_u,\tilde{\phi}_u)=\mathrm{csch\,}^2\theta
\end{equation}
Considering $e_2(\theta)=0$ from Proposition \ref{PropClassA} and \eqref{lemmaeq1}, we have 
$\partial_v (g_c(\tilde{\phi}_u,\tilde{\phi}_u))=0$.
Thus, we get the desired result.  
\end{proof}

The proof of the following corollary directly follows from the proof of Lemma \ref{space-likeSurfCanonicalParamLemma2}.
\begin{corry}
A space-like surface in $L^4_1(f,c)$ satisfies $A_\eta T=\lambda T$ for a smooth function $\lambda$ 
and a normal vector $\eta$ if and only if it can be locally parameterized by \eqref{space-likeSurfCanonicalParam} 
for an immersion $\tilde\phi$  satisfying \eqref{space-likeSurfCanonicalParamLemma2Eq1}. 
\end{corry}

\section{Class $\mathcal A$ surfaces in $L^4_1(f,0)$} 
In this section, we will give the local classification theorem of space-like class $\mathcal A$ surfaces for $c=0$.
Thus, we focus on the surface $M$ given in Lemma \ref{space-likeSurfCanonicalParamLemma2}.

\subsection{Local Classification Theorem}
Let  $\bar M$ denote the surface in $\mathbb E^3$ parameterized by $\tilde\phi$ with a unit normal vector field $\tilde N$ 
and $\tilde E,\tilde G$ stand for the coefficients of the first fundamental form of $\tilde M$, i.e.,
$\tilde E=g_0(\tilde\phi_u,\tilde\phi_u)$ and $\tilde G=g_0(\tilde\phi_v,\tilde\phi_v)$.

First, we  construct a geometrical frame field on $M$ as follows.
\begin{lemma}
Let $M$ be a space-like surface in $L^4_1(f,0)$ parameterized by \eqref{space-likeSurfCanonicalParam}. Then, 
there exists an orthonormal frame field $\{e_1, e_2, e_3, e_4\}$ on $M$ given by 
\begin{align}\label{space-likeSurfFrameFieldL41f0}
\begin{split}
e_1=&\frac1{\sqrt{-1+f^2\tilde E}}\partial_u,\quad   e_2=\frac1{f\sqrt{\tilde G}} \partial_v,\\
e_3=& \frac{1}{\sqrt{f^2 \tilde E^2-\tilde E}}\left(f \tilde E,\frac{1}{f }\tilde \phi_u\right),\\
e_4=& \frac1f (0,\tilde N).
\end{split}
\end{align}
Moreover, the functions $\theta$ and $h^4_{12}$ satisfy
\begin{eqnarray}
\label{space-likeSurfL41f0DecompTheta}\theta&=&-\sinh^{-1}\left(\frac1{\sqrt{-1+f^2\tilde E}}\right),\\
\label{connformspace-likeSurfL41f02}
h^4_{12}&=&\frac{f}{\sqrt{EG}} g_0(\tilde\phi_{uv},\tilde N).
\end{eqnarray}
\end{lemma}

\begin{proof}
By considering the parametrization of the immersion $\phi$ in \eqref{space-likeSurfCanonicalParam} and 
the first condition in \eqref{space-likeSurfCanonicalParamLemma2Eq1}, 
one can observe that the vector fields $e_1,e_2,e_3,e_4$ are as given in \eqref{space-likeSurfFrameFieldL41f0}. From Lemma \ref{Lemmafptsubman} and \eqref{space-likeSurfFrameFieldL41f0} we have $f'(t)=f'(u)$ whenever $(t,\tilde p)=\phi(u,v)\in M$.
Therefore, by a direct computation using \eqref{RWSTconn} and \eqref{space-likeSurfFrameFieldL41f0} we get desired results.
\end{proof}

Now, we consider the following examples of space-like class $\mathcal A$ surfaces in $L^4_1(f,0)$. 
\begin{example}\label{ClassAspace-likeExample1}
For some smooth functions $x_1,x_2$, consider the following space-like surface in $L^4_1(f,0)$ 
\begin{equation}\label{ClassAspace-likeExample1Param}
    \phi(u,v)=(u,x_1(u),x_2(u),v)
\end{equation}
with $-1+f^2(x_1'^2(u)+x_2'^2(u))>0$. Say $V(u)=\sqrt{x_1'^2(u)+x_2'^2(u)}$. 
Then, the equation \eqref{space-likeSurfFrameFieldL41f0} turns into
\begin{align}\label{MinClassACase1Eq1}
\begin{split}
e_1=&\frac1{\sqrt{-1+f^2V^2}}\partial_u,\quad   e_2=\frac1{f} \partial_v,\\
e_3=& \frac{1}{f V\sqrt{f^2 V^2-1}}\left(f^2 V^2,{x_1'},{x_2'},0\right),\\
e_4=& \frac1{fV} (0,-x_2',x_1',0)
\end{split}
\end{align}
and $\displaystyle{\sinh\theta=-\frac{1}{\sqrt{-1+f^2V^2}}}$. Thus, it can be seen that $e_2(\theta)=0$ and $\nabla^\perp_{e_2}e_4=0$.
From Proposition \ref{PropClassA}, we get that $M$ is a space-like class $\mathcal A$ surface in $L^4_1(f,0)$.
For later use, we also give the nonzero components of the second fundamental form $h$ of $M$ as
\begin{align}\label{MinClassACase1Eq2}
\begin{split}
h^3_{11}=&\frac{f V'-V \left(f^2 V^2-2\right) f'}{\left(f^2 V^2-1\right)^{3/2}},\qquad
h^3_{22}=-\frac{V f'}{\sqrt{f^2 V^2-1}},\\
h^4_{11}=&\frac{f \left(x_2' x_1''-x_1' x_2''\right)}{V-f^2 V^3},\qquad
h^4_{22}=0.  
\end{split}
\end{align}
\end{example}

\begin{example}\label{ClassAspace-likeExample2}
Let $\alpha:I_v\to\mathbb S^2$ be an arc-length parameterized curve with unit normal $n$ and curvature $\kappa$.
Consider the following space-like surface in $L^4_1(f,0)$ given by
\begin{equation}\label{ClassAspace-likeExample2Param}
    \phi(u,v)=(u,\phi_1(v)\alpha(v)+\phi_2(u,v)\alpha'(v)+\phi_3(u,v)n(v))
\end{equation}
where $\phi_1(v)$, $\phi_2(u,v)$ and $\phi_3(u,v)$ are smooth functions defined by
\begin{align}\label{ClassAspace-likeExample2ParamA2A3}
\begin{split}
       \phi_2(u,v)=&\int_{u_0}^u R(\xi)\sin(\tau(\xi,v))d\xi +\psi_1(v)\\
    \phi_3(u,v)=&\int_{u_0}^u R(\xi)\cos(\tau(\xi,v))d\xi+\psi_2(v) 
\end{split}
\end{align}
for some smooth functions $R$, $\tau$, $\psi_1$ and $\psi_2$ satisfying $\tau_v=\kappa$ and 
\begin{equation}\label{ClassAspace-likeExample2ParamA4A5}
\psi_1'=\kappa \psi_2-\phi_1,\qquad \psi_2'=-\kappa \psi_1.
\end{equation}
Since it is space-like, we have $(-1+f^2R^2)(\phi_1'-\phi_2)>0$. 
By a direct computation, one can obtain that 
\begin{eqnarray}
  e_2=\frac 1{f(\phi_1'-\phi_2)}\partial_v&\quad&  e_4=\frac 1{f}(0,\cos\tau \alpha'-\sin\tau n),\\
\label{ClassAspace-likeExample2theta}   \sinh\theta=-\frac{1}{\sqrt{-1+f^2R^2}}.&&
\end{eqnarray}
It can be seen that $e_2(\theta)=0$ and $\nabla^\perp_{e_2}e_4=0$.
Thus, Proposition \eqref{PropClassA} implies that $M$ is a space-like class $\mathcal A$ surface in $L^4_1(f,0)$. 
\end{example}

Moreover, if $\kappa=0$, then the surface described in 
Example \ref{ClassAspace-likeExample2} turns into the surface of revolution given by the following example.

\begin{example}\label{minimalclassAimmersionREMARK}
Consider the following space-like surface parameterized by
\begin{equation} 
\label{minimalclassAimmersion}
\phi(u,v)=(u,\zeta_1(u)\cos{v}, \zeta_1(u)\sin{v}, \zeta_2(u))
\end{equation}
for some smooth functions $\zeta_1>0$ and $\zeta_2$. 
Let define $V(u)=\sqrt{\zeta _1'^2(u)+\zeta _2'^2(u)}$.
Then, \eqref{space-likeSurfFrameFieldL41f0} turns into
\begin{align*}
\begin{split}
e_1=&\frac1{\sqrt{-1+f^2 V^2}}\partial_u,\quad   e_2=\frac1{f\zeta_1} \partial_v,\\
e_3=& \frac{1}{f V \sqrt{-1+f^2 V^2}}\left(f^2 V^2 ,\zeta _1'\cos  v ,\zeta _1'\sin  v ,\zeta _2'\right),\\
e_4=& \frac1{fV} (0,\zeta _2'\cos  v ,\zeta _2'\sin  v ,-\zeta _1'),
\end{split}
\end{align*}
By a direct computation, we obtain the nonzero coefficients of the second fundamental form $h$ of $M$ as
\begin{align}\label{MinClassACase2P2Eq2}
\begin{split}
h^3_{11}=&\frac{V  \left(2-f ^2 V ^2\right) f' +f  V' }{\left(-1+f ^2 V ^2\right)^{3/2}},\qquad
h^3_{22}=\frac{V ^2 \sqrt{-1+f ^2 V ^2} f' +\zeta _1' }{V -f ^2 V ^3},\\
h^4_{11}=&\frac{f  \left(\zeta _2'  \zeta _1'' -\zeta _1'  \zeta _2'' \right)}{V  \left(-1+f ^2 V ^2\right)},\qquad
h^4_{22}=-\frac{\zeta _2' }{V  \sqrt{-1+f ^2 V ^2}}.  
\end{split}
\end{align}
\end{example}

Now, we give the local classification of space-like class $\mathcal A$ surfaces in $L^4_1(f,0)$.

\begin{theorem}
\label{theoremclassAf0}
A space-like surface in $L^4_1(f,0)$ is a class $\mathcal A$ surface if and only if it is locally congruent to one of the following surfaces:
\begin{itemize}
    \item[(i)] The cylinder described in Example \ref{ClassAspace-likeExample1},
    \item[(ii)] The surface described in Example \ref{ClassAspace-likeExample2}.
\end{itemize}
\end{theorem}

\begin{proof} 
Assume that $M$ is a space-like class $\mathcal A$ surface in $L^4_1(f,0)$ and let $p\in M$. 
Then, Lemma \ref{space-likeSurfCanonicalParamLemma2} implies that 
there exists local coordinates $(u,v)$ on a neighborhood $\mathcal N_p$ of $p$ 
which can be parameterized by \eqref{space-likeSurfCanonicalParam} immersion $\tilde\phi$ satisfying \eqref{space-likeSurfCanonicalParamLemma2Eq1}. 
Then, 
from Proposition \ref{PropClassA}, we have $\tilde g(h(e_1,e_2),e_4)=h^4_{12}=0$.
Thus, the equation \eqref{connformspace-likeSurfL41f02} gives 
$g_0(\tilde\phi_{uv},\tilde N)=0$. 
Therefore, we express 
$$\tilde\phi_{uv}=\frac{g_0(\tilde\phi_{uv},\tilde\phi_{v})}{g_0(\tilde\phi_{v},\tilde\phi_{v})}\tilde\phi_{v}=\frac{\tilde G_u}{2\tilde{G}}\tilde\phi_{v}$$
from which we get 
\begin{equation}\label{SpClassAL41fNSCond}
  \tilde{\phi}_v=\sqrt{\tilde{G}(u,v)}\alpha (v)  
\end{equation}
for an $\mathbb R^3$-valued function $\alpha$ satisfying $g_0( \alpha, \mathbf \alpha)=1$. 
Then, there occur two following cases according to the function $\alpha$.

\textit{Case (i.)} $\alpha$ is a constant vector in $\mathbb R^3$. 
In this case, up to a suitable isometry, we may assume $\alpha=(0,0,1)$. 
Considering \eqref{SpClassAL41fNSCond}, the case (i) of the theorem can be obtained. 

\textit{Case (ii.)} $\alpha$ is a non constant vector in $\mathbb R^3$, that is, $\alpha'\neq 0$. 
In this case, by re-defining $v$ properly, one may assume $g_0( \alpha', \alpha')=1$. 
Thus, $\alpha:I_v\to\mathbb S^2$ is an arc-length parameterized spherical curve. 
Let $\kappa$ and $n$ be a curvature and a normal vector of $\alpha$, 
where we have $\alpha''=\kappa n-\alpha$ and $n'=-\kappa \alpha'$. Then, the immersion $\tilde\phi$ in $\mathbb{E}^3$ can be written as
\begin{equation}\label{SpClassAL41fNSCase2Eq1}
\tilde\phi(u,v)= \phi_1(u,v)\alpha(v) +\phi_2(u,v)\alpha'(v)+\phi_3(u,v)n(v)
\end{equation}
for some smooth functions $\phi_1,\phi_2,\phi_3$. 
Since $\tilde{g}(\tilde\phi_u,\tilde\phi_v)=0$, 
the equations \eqref{SpClassAL41fNSCond} and 
\eqref{SpClassAL41fNSCase2Eq1} give 
$\displaystyle{\frac{\partial \phi_1}{\partial u}=0}$ which implies $\phi_1=\phi_1(v)$.
Considering
\eqref{space-likeSurfL41f0DecompTheta} and $e_2(\theta)=0$ together, we have $\tilde E(u,v)=R^2(u)$ 
for a non-vanishing function $R$. Thus,
we get
\begin{equation}\label{SpClassAL41fNSCase2Eq2}
\frac{\partial \phi_2}{\partial u}=R\sin\tau,\quad 
\frac{\partial \phi_3}{\partial u}=R\cos\tau,
\end{equation}
for a smooth function $\tau=\tau(u,v)$. 
Integrating \eqref{SpClassAL41fNSCase2Eq2}, we get the equations \eqref{ClassAspace-likeExample2ParamA2A3} 
for $\phi_2$ and $\phi_3$
with some functions $\psi_1$ and $\psi_2$. 
On the other hand, by a direct computation using \eqref{ClassAspace-likeExample2Param}, 
\eqref{SpClassAL41fNSCase2Eq2} and $g_0(\tilde\phi_u,\tilde\phi_v)=0$, we obtain
\begin{equation}\label{SpClassAL41fNSCase2Eq3}
\frac{\partial \phi_2}{\partial v}+\phi_1-\kappa \phi_3=0,\quad 
\frac{\partial \phi_3}{\partial v}+\kappa \phi_2=0.
\end{equation}
Taking the derivative of \eqref{SpClassAL41fNSCase2Eq2} and \eqref{SpClassAL41fNSCase2Eq3} with respect to $v$ and $u$, 
we get $\frac{\partial \tau}{\partial v}=\kappa$  and \eqref{ClassAspace-likeExample2ParamA2A3} for some smooth functions $\psi_1$ and $\psi_2$. 
By using \eqref{ClassAspace-likeExample2ParamA2A3} and \eqref{SpClassAL41fNSCase2Eq3}, we obtain \eqref{ClassAspace-likeExample2ParamA4A5}.
Thus, $M$ is congruent to the surface given in case (ii) of the theorem.
Hence, the proof of the necessary condition is completed. 

The proof of the sufficient condition is obtained from Example \ref{ClassAspace-likeExample1} and Example \ref{ClassAspace-likeExample2}. 
\end{proof}

\subsection{Minimal Surfaces}
In this subsection, as an application of Theorem \ref{theoremclassAf0} we study minimal class $\mathcal A$ surfaces in $ L^4_1(f,0)$.
First, we will focus the surface given in (i) of Theorem \ref{theoremclassAf0}.
\begin{proposition}\label{classminimalCase1}
The class $\mathcal A$ space-like surface described in (i) of Theorem \ref{theoremclassAf0} is minimal 
if and only if it is congruent to the surface given by 
\begin{equation}\label{MinClassACase1Eq4}
    \phi(u,v)=(u,x_1(u),c_1x_1(u)+c_2,v)
\end{equation}
for constants $c_1,\ c_2$ and a function $x_1$ satisfying 
\begin{equation}\label{MinClassACase1Eq3}
  x_1=\int_{u_0}^u\frac{d\xi}{f(\xi)\sqrt{c_3 f^4(\xi)+(c_1^2+1)}}
\end{equation} 
where $u_0, c_3$ are real constants.
\end{proposition}
\begin{proof}
Assume that $M$ is a minimal class $\mathcal A$ surface described in (i) of Theorem \ref{theoremclassAf0},
that is, $ h^3_{11}+ h^3_{22}=h^4_{11}+ h^4_{22}=0$. 
From the equations in \eqref{MinClassACase1Eq2}, we find $x_2'x_1''-x_1'x_2''=0$ whose solution is given by  
$x_2(u)=c_1x_1(u)+c_2$ for some constants $c_1,c_2$. 
By a further computation, considering $ h^3_{11}+ h^3_{22}=0$ and \eqref{MinClassACase1Eq2} 
we obtain the following differential equation
$$fx_1''-2f^2(1+c_1^2)f'x_1'{}^3 +3f'x_1'=0.$$
Then, the solution of this equation is given by \eqref{MinClassACase1Eq3}. 
\end{proof}
\begin{rem}\label{classminimalCase1Rem}
The surface given by \eqref{MinClassACase1Eq4} lies on a totally geodesic hypersurface of $L^4_1(f,0)$.
\end{rem}    

Now, we will examine the surface given in (ii) of Theorem \ref{theoremclassAf0} with minimality condition. 
\begin{lemma}
\label{lemclassminimalCase2}
The class $\mathcal A$ space-like surface described in (ii) of Theorem \ref{theoremclassAf0} 
is a minimal surface if and only if $\kappa=0$. 
\end{lemma}
\begin{proof}
Assume that $M$ is a minimal class $\mathcal A$ surface described in (ii) of Theorem \ref{theoremclassAf0}. 
Let $\bar M$ denote the surface in $\mathbb E^3$ parameterized by 
$$\tilde\phi(u,v)=\phi_1(v)\alpha(v)+\phi_2(u,v)\alpha'(v)+\phi_3(u,v)n(v).$$
By a direct computation, we obtain that the principal directions of $\bar M$ are 
$\bar e_1=\frac 1R\tilde\phi_*(\partial_u)$ and $\bar e_2=\frac 1{\phi_1'-\phi_2}\tilde\phi_*(\partial_v)$ 
with the corresponding principal curvatures $k_1,\ k_2$ given by
\begin{equation}
\label{MinClassACase2Eq2}
k_1=\frac{\tau_u}{R},\qquad k_2=\frac{\cos\tau}{\phi_1'-\phi_2}.
\end{equation}
Since $\tau_v=\kappa(v)$, we get $\bar{e_2}(k_1)=0$ which means $k_1=k_1(u)$. 
By a direct computation, we obtain that $k_1$ and $k_2$ are related with
\begin{equation}\label{MinClassACase2Eq3}
\bar e_1(k_2)=\frac{\sin\tau}{\phi_1'-\phi_2}(k_2-k_1).
\end{equation}
On the other hand, \eqref{space-likeSurfFrameFieldL41f0} implies
\begin{equation}\label{MinClassACase2Eq1}
   e_4= \frac1{f} (0,\cos\tau \alpha'-\sin\tau n). 
\end{equation}
Considering \eqref{MinClassACase2Eq2} and \eqref{MinClassACase2Eq1},
we obtain
\begin{equation}
\label{MinClassACase2Eq4}
h^4_{11}+h^4_{22}=\frac{fR^2k_1}{-1+f^2R^2}+\frac{k_2}{f}.
\end{equation}
Since $M$ is a minimal surface, that is, $h^4_{11}+h^4_{22}=0$, 
$\displaystyle{k_2=\frac{f^2R^2k_1}{1-f^2R^2}}$ which implies $k_2=k_2(u)$.
Considering this with \eqref{MinClassACase2Eq2} and \eqref{MinClassACase2Eq3}, we obtain $\tau_v=\kappa=0$. 
\end{proof}

Hence, we give the following proposition:
\begin{proposition}
\label{classminimalCase2}
The class $\mathcal A$ surface described in (ii) of Theorem \ref{theoremclassAf0} is minimal if and only if it is congruent to the surface given by 
\eqref{minimalclassAimmersion} in Example \ref{minimalclassAimmersionREMARK}
for some smooth functions $\zeta_1$ and $\zeta_2$ satisfying 
\begin{align}\label{classminimalCase2Eq}
\begin{split}
f\zeta _1''=& {f' \zeta _1' \left(2 f^2 \left(\zeta _1'{}^2+\zeta _2'{}^2\right)-3\right)+\sqrt{f^2 \left(\zeta _1'{}^2+\zeta _2'{}^2\right)-1}},\\
f\zeta _2''=& {f' \zeta _2' \left(2 f^2 \left(\zeta _1'{}^2+\zeta _2'{}^2\right)-3\right)}.
\end{split}
\end{align}
\end{proposition}

\begin{proof}
Suppose that $M$ is a class $\mathcal A$ surface described in (ii) of Theorem \ref{theoremclassAf0} and it is minimal.
Thus, \eqref{lemclassminimalCase2} implies that $M$ is congruent to the surface parameterized by \eqref{minimalclassAimmersion} 
for some smooth functions $\zeta_1>0$ and $\zeta_2$.
By considering \eqref{MinClassACase2P2Eq2}, we obtain the equations in \eqref{classminimalCase2Eq}.
\end{proof}

By combining Proposition \ref{classminimalCase1}, Remark \ref{classminimalCase1Rem} and Proposition \ref{classminimalCase2}, 
we obtain the following result:
\begin{theorem}
\label{classminimalTHM}
Let $M$ be a space-like surface in $L^4_1(f,0)$ which has no open part lying on a totally geodesic hypersurface of $L^4_1(f,0)$. 
Then, $M$ is minimal and class $\mathcal A$ surface if and only if it is locally congruent to the surface \eqref{minimalclassAimmersion}
for some smooth functions $\zeta_1$ and $\zeta_2$ satisfying \eqref{classminimalCase2Eq}.
\end{theorem}
  
\section{space-like Surfaces in $L^4_1(f,0)$ satisfying $\nabla^\perp \eta=0$}
In this section, we consider space-like surfaces in $L^4_1(f,0)$ having a parallel normal vector field $\eta$ 
in \eqref{RWTS-etaandTDef}. 
First, we will give the following lemma to obtain the classification theorem.
\begin{lemma}
\label{LMEtaParallel}
Let $M$ be a space-like surface in $L^4_1(f,c)$. 
Then, the vector field $\eta$ in \eqref{RWTS-etaandTDef} is parallel if and only if $M$ is a class $\mathcal A$ surface 
in $L^4_1(f,c)$.
\end{lemma}

\begin{proof}
Suppose that $M$ is a space-like surface in $L^4_1(f,c)$ and $\nabla^\perp_X\eta=0$ for any tangent vector field $X$. 
From the equation \eqref{RWTS-etaandTDefNew}, we know that $\eta=\cosh{\theta}e_3$. 
Thus, we get that $\theta$ is a nonzero constant and $\nabla^\perp_Xe_3=0$. 
Note that when $\theta$ is zero, $\frac{\partial}{\partial t}=e_3$ which means $T=0$. We omit this case.
From Proposition \ref{PropClassA}, it can seen that $M$ is a class $\mathcal A$ surface in $L^4_1(f,c)$. 
\end{proof}

By using Lemma \ref{LMEtaParallel}, we will state the following Theorem. 

\begin{theorem}\label{THMEtaParallel}
A space-like surface in $L^4_1(f,0)$ having the parallel vector field $\eta$ in \eqref{RWTS-etaandTDef} is locally congruent 
to one of the following surfaces:
\begin{itemize}
    \item[(i)] The cylinder given by 
    \begin{equation}\label{THMEtaParallelC1}
    \phi(u,v)=\left(u,c_1\int_{u_0}^u\frac{d\xi}{f(\xi)},c_2\int_{u_0}^u\frac{d\xi}{f(\xi)}+c_3,v\right)
\end{equation}
for some constants $u_0,\ c_1,\ c_2$ and $c_3$.
    \item[(ii)] The surface parameterized by \eqref{ClassAspace-likeExample2Param}, where $A_1$ is a smooth function and the functions
    $A_2,A_3$ are given by
\begin{align}\label{THMEtaParallelC2}
\begin{split}
    A_2(u,v)=&c\sin\tau(v)\int_{u_0}^u \frac{d\xi}{f(\xi)} +A_4(v)\\
    A_3(u,v)=&c\cos\tau(v)\int_{u_0}^u \frac{d\xi}{f(\xi)}+A_5(v) 
\end{split}
\end{align}
for constants $c,u_0$, smooth functions $\tau, A_4$ and $A_5$ satisfying
\eqref{ClassAspace-likeExample2ParamA4A5} and $\tau'(v)=\kappa(v)$.
\end{itemize}
\end{theorem}

\begin{proof}
Suppose that $M$ is a space-like surface in $L^4_1(f,0)$ with parallel vector field $\eta$.
Then, Lemma \eqref{LMEtaParallel} implies that $M$ is a class $\mathcal A$ surface.
Thus, we have two surfaces given in Theorem \ref{theoremclassAf0}.
Now, we are going to study these surfaces, separately. 

    \textit{Case (i)}. $M$ is congruent to the cylinder described in Example \ref{ClassAspace-likeExample1} for some smooth functions $x_1,\ x_2$. 
    In this case, since $\theta$ is constant, we have
    \begin{equation}\label{THMEtaParallelC1Eq1}
    x_1'{}^2(u)+x_2'{}^2(u)=\frac{c_1^2}{f^2(u)}
    \end{equation}
    for a constant $c_1>0$. From the equation \eqref{LemmaCurvProps6}, we get $h^4_{11}=0$. 
    Thus, the equation \eqref{MinClassACase1Eq2} gives $x_2' x_1''-x_1' x_2'=0$ 
    which yields $x_2(u)=a_1x_1(u)+a_2$ for some constants $a_1,a_2$.
    By combining this equation with \eqref{THMEtaParallelC1Eq1}, we get $x_1'=\frac{c_1}{f}$.
    Then, \eqref{ClassAspace-likeExample1Param} turns into \eqref{THMEtaParallelC1}. Hence, we have the case (i) of the theorem.
    
    \textit{Case (ii)}.  The surface described in Example \ref{ClassAspace-likeExample2}.  In this case, similar to the case (i), we have 
    \begin{equation}\label{THMEtaParallelC2Eq1}
        R=\frac{c}{f}.
    \end{equation}
Using \eqref{MinClassACase2Eq1} and $h^4_{11}=0$, we get $\tau_u=0$. 
Thus, \eqref{ClassAspace-likeExample2ParamA2A3} turns into \eqref{THMEtaParallelC2}.
Hence, we have case (ii) of the theorem. 
  
The converse can be shown through a straightforward calculation.
\end{proof}

\section*{Acknowledgements}
The authors declare that they have no conflict of interest.
This work was carried out during the 1001 project supported by the Scientific 
and Technological Research Council of T\"urkiye (T\"UB\.ITAK)  (Project Number: 121F352).

\end{document}